\documentclass[a4paper]{article}
\usepackage{amssymb,amsmath,amsthm,verbatim,IEEEtrantools,macros,url,authblk}
\begin{document}
\title{Three-term arithmetic progressions in subsets of $\mathbb{F}_q^{\infty}$ of large Fourier dimension}
\author{Robert Fraser}
\affil{University of Edinburgh}
\maketitle
\begin{abstract}
We show that subsets of $\mathbb{F}_q^{\infty}$ of large Fourier dimension must contain three-term arithmetic progressions. This contrasts with a construction of Shmerkin of a subset of $\mathbb{R}$ of Fourier dimension $1$ with no three-term arithmetic progressions.
\end{abstract}
\section{Introduction}
In a recent paper, Ellenberg and Gijswijt \cite{EllenbergGijswijt17} have shown that, for any odd prime $q$, there exists $r < q$ such that a subset of $\mathbb{F}_q^d$ with at least $r^d$ elements must contain a three-term arithmetic progression. This means that, in contrast to the case for finite cyclic groups, where Behrend \cite{Behrend46} constructed a counter example, a sparse subset of a finite vector space of sufficiently large (but not full) dimension must contain a three-term arithmetic progression.

In this note, we will consider what happens in the vector space $\mathbb{F}_q^{\infty}$- a vector space of infinite dimension over $\mathbb{F}_q^d$. By $\mathbb{F}_q^{\infty}$, we mean the vector space consisting of infinite sequences of elements of $\mathbb{F}_q^{d}$ with the product topology. This is a compact abelian group that is isomorphic to the additive group of $\mathbb{F}_q[[t]]$, the ring of formal power series over $\mathbb{F}_q$. 

In light of the result of Ellenberg and Gijswijt, one may be tempted to guess that a subset of $\mathbb{F}_q^{\infty}$ of full Hausdorff dimension must contain a three-term arithmetic progression; however, this has been shown not to be the case \cite{Fraser18}. The construction was inspired by a similar construction of Keleti \cite{Keleti98} of a subset of $\mathbb{R}$ of Hausdorff dimension that does not contain any solutions to $x_4 - x_3 = x_2 - x_1$ with $x_1 \neq x_2$ and $x_3 \neq x_4$. Because size in the sense of Hausdorff dimension is not enough to guarantee the existence of a three-term arithmetic progression, some additional condition, such as a Fourier decay condition, is needed.

In the real-variable setting, \L aba and Pramanik \cite{LabaPramanik09} have shown that a subset of $\mathbb{R}$ supporting a measure satisfying a Fourier decay condition as well as a ball condition depending on the rate of Fourier decay must contain a three-term arithmetic progression. However, Shmerkin \cite{Shmerkin17} has constructed a subset of $\mathbb{R}$ of Fourier dimension $1$ not containing any three-term arithmetic progressions. Shmerkin's construction relied on the Behrend example \cite{Behrend46} of a large subset of $\{1, 2, \ldots, N\}$ that does not contain a three-term arithmetic progression. Because the result \cite{EllenbergGijswijt17} of Ellenberg and Gijswijt implies that no such example can exist for finite vector spaces, it seems sensible to guess that a subset of $\mathbb{F}_q^{\infty}$ with large Fourier dimension must contain a three-term arithmetic progression. This is exactly what we will show:

\begin{mythm}\label{mainthm}
Let $q$ be an odd prime. Then for any $1 > \beta > 2/3$, there exists $\alpha < 1$ depending only on $\beta$ and $q$ with the following property: let $E$ be a compact subset of $\mathbb{F}_q^{\infty}$ supporting a probability measure $\mu$ such that for some positive constants $C_1$ and $C_2$:
\begin{enumerate}
\item There exists $E' \subset E$ such that $\mu(E') > 0$ and for all balls $B \subset \mathbb{F}_q^{\infty}$,
\[\mu'(B) \leq C_1 \text{rad}(B)^{\alpha}.\]
Here, $\mu'$ is the measure $\mu$ restricted to $E'$.
\item \[\hat \mu(\xi) \leq C_2 |\xi|^{- \beta/2}\]
\end{enumerate}
Then the set $E$ must contain a three-term arithmetic progression. The second condition implies the first for all $\alpha < \beta$. If $\beta$ is sufficiently close to $1$ depending on $q$, then the first condition is unnecessary.
\begin{IEEEeqnarray*}{rCl}
\end{IEEEeqnarray*}
\end{mythm}
This differs from the \L aba-Pramanik result \cite{LabaPramanik09} because the value $\alpha$ does not depend on the constants $C_1$ and $C_2$. This allows us to drop the first assumption provided that $\beta$ is sufficiently close to $1$. The counterexample of Shmerkin \cite{Shmerkin17} shows that this assumption cannot be dropped in the Euclidean setting.

In order to properly interpret this theorem, we need to discuss some of the basic properties of the Fourier transform on $\mathbb{F}_q^{\infty}.$
\subsection{Acknowledgements}
The author would like to thank Jonathan Hickman for his encouragement and advice. The author would also like to thank Cosmin Pohoata for pointing out the history behind Lemma \ref{Varnavides}.

This material is based on work supported by the NSF under Award No. 1803086.
\section{Fourier Analysis on $\mathbb{F}_q^{\infty}$}
\subsection{The abelian groups $\mathbb{F}_q^{\infty}$ and $\widehat{\mathbb{F}_q^{\infty}}$}
Much of the material in this section can be found in Taibleson's book \cite{Taibleson75}. Let $q$ be an odd prime, and let $\mathbb{F}_q^{\infty}$ be the group
\[\prod_{j=1}^{\infty} \mathbb{F}_q\]
equipped with the product topology. With respect to this topology, $\mathbb{F}_q^{\infty}$ is a compact abelian group. The topology on $\mathbb{F}_q^{\infty}$ is induced by an absolute value: given an element $x = (x_0, x_1, x_2, \ldots)$ of $\mathbb{F}_q^{\infty}$, we define $|x| = q^{-j}$, where $j$ is the index of the first nonzero component of $x$. If $x = 0$, then we take $|x| = 0$. 
There is a natural projection $\pi_d : \mathbb{F}_q^{\infty} \to \mathbb{F}_q^d$ given by $\pi_d(x) = (x_0, \ldots, x_{d-1})$. Note that for any $d^* > d$, there is a natural projection from $\mathbb{F}_q^{d^*} \to \mathbb{F}_q^d$; we will abuse notation and also use $\pi_d$ for this projection. As for $\mathbb{F}_q^{d}$, we define an absolute value on $\mathbb{F}_q^{d}$ by $|(x_0, \ldots, x_{d-1})| = q^{-j}$, where $j$ is the index of the first nonzero component of $(x_0, \ldots, x_{d-1})$, and $|(0,0,\ldots,0)| = 0$. Notice that if $x \in \mathbb{F}_q^{\infty}$ is such that $|\pi_d(x)| > 0$, then $|\pi_d(x)| = |x|$. 

The compact abelian group $\mathbb{F}_q^{\infty}$ is equipped with a Haar probability measure $dx$. This measure assigns a measure of $q^{-j}$ to any closed ball of radius $q^{-j}$. The pushforward of this measure under $\pi_d$ yields the uniform probability measure on $\mathbb{F}_q^d$.

The Fourier character group $\widehat{\mathbb{F}_q^d}$ of $\mathbb{F}_q^d$ is isomorphic to $\mathbb{F}_q^d$ as an abelian group. We will write $(\xi_1, \ldots, \xi_{d})$ for a typical character on $\mathbb{F}_q^d$ (notice that the indexing will start from $1$ instead of $0$). We define an absolute value on $\widehat{\mathbb{F}_q^d}$ by $|\xi| = q^j$, where $j$ is the maximum index of a nonzero component of $(\xi_1, \ldots, \xi_{d})$. The product $(\xi_1, \ldots, \xi_{d}) \cdot (x_0, \ldots, x_{d-1})$ is defined by $\xi_1 x_0 + \cdots + \xi_{d} x_{d-1}$, which is defined as an element of $\mathbb{F}_q$. We can therefore make sense of $\exp(\frac{2 \pi i }{q} \xi \cdot x)$, which will be written as 
\begin{equation}\label{eqdef}
e_q(\xi \cdot x) = \exp(\frac{2 \pi i}{q} \xi \cdot x).
\end{equation}
 This describes the action of $\widehat{\mathbb{F}_q^d}$ on $\mathbb{F}_q^d$. 

The Fourier character group $\widehat{\mathbb{F}_q^{\infty}}$ consists of sequences of the form $\xi = (\xi_1, \xi_2, \xi_3, \ldots)$ where only finitely many $\xi_j$ are nonzero. The absolute value $|\xi|$ of $\xi$ is given by $q^j$, where $j$ is the largest index of a nonzero component of $\xi$, with $|0|$ taken to be $0$. Because all of the components of $\xi$ after the $j$th component are zero, we can define a product $\xi \cdot x$ for $\xi \in \widehat{\mathbb{F}_q^{\infty}}$ and $x \in \mathbb{F}_q^{\infty}$ as the finite sum 
\[\sum_{k=1}^j \xi_k x_{k-1}\]
which makes sense as an element of $\mathbb{F}_q$. We can thus define $e_q(\xi \cdot x)$ as before, giving the action of $\widehat{\mathbb{F}_q^{\infty}}$ on $\mathbb{F}_q^{\infty}$. Notice that each element of $\widehat{\mathbb{F}_q^{\infty}}$ can be viewed as an element of $\widehat{\mathbb{F}_q^d}$ where $d \geq j$ and $q^j = |\xi|$. In this sense, every element of $\widehat{\mathbb{F}_q^{\infty}}$ can be viewed as an element of $\widehat{\mathbb{F}_q^j}$ for some finite $j$. In fact, if $|\xi| \leq q^d$, then $\xi \cdot x = \xi \cdot \pi_d(x)$. In other words, the function $x \mapsto \xi \cdot x$ is constant on closed balls of radius $|\xi|^{-1}$ for $\xi \neq 0$.

Given $d^* > d$, and $x \in \mathbb{F}_q^{d^*}$, we can write $x = (x_0, x_1, \ldots, x_{d^* - 1})$ as a sum $x = x' + x''$, where 
\begin{IEEEeqnarray*}{rCl}
x' & = & (x_0, \ldots, x_{d - 1}, 0, \ldots, 0) \\  
x''& = & (0, \ldots, 0, x_d, \ldots, x_{d^*-1}).
\end{IEEEeqnarray*} 
We call this the \textbf{order $d$ decomposition} on $\mathbb{F}_q^{d^*}$. Similarly, given $\xi \in \widehat{\mathbb{F}_q^{d^*}}$, we can write $\xi = \xi' + \xi''$, where 
\begin{IEEEeqnarray*}{rCl}
\xi' & = & (\xi_1, \ldots, \xi_d, 0, \ldots, 0) \\
\xi'' & = & (0, \ldots, 0, \xi_{d+1}, \ldots, \xi_{d^*}).
\end{IEEEeqnarray*} 
We will call this the \textbf{order $d$ decomposition} of $\xi$. We note some trivial facts about these order $d$ decompositions. First, we observe that  $|x''| \leq q^{-d}$ and $|\xi'| \leq q^d$. We have $|\xi''| \geq q^{d+1}$ unless $\xi'' = 0$. We also have that $(x' + x'') \cdot (\xi' + \xi'') = x' \cdot \xi' + x'' \cdot \xi''$. 
\subsection{The Fourier transform on $\mathbb{F}_q^{\infty}$}

The Fourier transform of an $L^1$ function $f : \mathbb{F}_q^{\infty} \to \mathbb{C}$ is given by
\[\hat f(\xi) = \int f(x) e_q (x \cdot \xi) \, dx\]
where $dx$ is the Haar measure on $\mathbb{F}_q^{\infty}$ and $e_q$ is as defined in \eqref{eqdef}. The Fourier transform of a finite measure $\mu$ on $\mathbb{F}_q^{\infty}$ is 
\[\hat \mu(\xi) = \int e_q(x \cdot \xi) \, d \mu(x)\]

The Fourier transform of a function $f : \mathbb{F}_q^d \to \mathbb{C}$ is 
\[\hat f(\xi) = \sum_{x \in \mathbb{F}_q^d} f(x) e_q(x \cdot \xi)\]
Notice that, if $\mu$ is a measure on $\mathbb{F}_q^{\infty}$ and $\mu_d$ is the pushforward of $\mu$ under $\pi_d$ (which can be interpreted as a function on $\mathbb{F}_q^d$), and $|\xi| \leq q^d$, then we have (by conflating $\xi \in \mathbb{F}_q^{\infty}$ with $\xi \in \mathbb{F}_q^d$ as above)
\begin{IEEEeqnarray*}{rCl}
\hat\mu_d(\xi) & = & \sum_{x \in \mathbb{F}_q^d} \mu_d(x) e(x \cdot \xi) \\
& = & \sum_{x_0 \in \mathbb{F}_q^d} \int_{\pi_d(x) = x_0} e(x_0 \cdot \xi)\, d \mu(x) \\
& = & \sum_{x_0 \in \mathbb{F}_q^d} \int_{\pi_d(x) = x_0} e(x \cdot \xi) \, d \mu(x) \\
& = & \hat \mu(\xi)
\end{IEEEeqnarray*}
This means that the Fourier coefficients of a measure $\mu(\xi)$ where $|\xi| \leq q^d$ can be computed directly in the finite vector space $\mathbb{F}_q^d$ without passing to the limit.
\subsection{Hausdorff and Fourier dimension of subsets of $\mathbb{F}_q^{\infty}$}
A good general reference for Hausdorff and Fourier dimensions in Euclidean spaces is \cite{Mattila15}. The notion of Fourier dimension occurring in this section is the $\mathbb{F}_q^{\infty}$ equivalent of the Euclidean Fourier dimension. Most of the material in this section appears in the thesis of Christos Papadimitropoulos \cite{Papadimitropoulos10C}.

Because $\mathbb{F}_q^{\infty}$ is a metric space, we can define the Hausdorff dimension of compact subsets of $\mathbb{F}_q^{\infty}$ in the usual manner. We will briefly review this definition now. 

For a compact set $E \subset \mathbb{F}_q^{\infty}$,  $t > 0$, define a $t$-covering of $E$ to be a covering of $E$ by closed balls of radius at most $t$. Define the $s$-dimensional $t$-Hausdorff content of $E$ as follows:
\[\mathcal{H}_t^d(E) :=  \inf_{\mathcal{B} \text{ $t$-covering of $E$}} \sum_{B \in \mathcal{B}} \text{rad}(B)^s\]

The value of $\mathcal{H}_t^d(E)$ increases as $t \to 0$ because the infimum is taken over a smaller family of coverings. We define
\[\mathcal{H}^s(E) := \sup_{t > 0} \mathcal{H}_t^s(E),\]
with the understanding that this supremum may be infinite. 

$H^s(E)$ is a non-increasing function of $s$. In fact, $H^s(E)$ will be equal to either $0$ or $\infty$ except for at most one value of $s$. Let $s_0 = \sup \{s : \mathcal{H}^s(E) = \infty\}$. Then $s_0$ is called the Hausdorff dimension of the set $E$. Note that $H^s(E)$ may be equal to $0$, $\infty$, or a finite non-zero value.

Frostman's Lemma relates the Hausdorff dimension of a compact subset $E$ of $\mathbb{F}_q^{\infty}$ to the ball condition of measures supported on the set $E$. In fact, this statement holds without the assumption that $E$ is compact, but that is all we will need.

The following version of Frostman's lemma can be found in Mattila \cite{Mattila95} as Theorem 8.17.

\begin{mylem}[Frostman's Lemma on Compact Metric Spaces]
Let $X$ be a compact metric space such that $\mathcal{H}^s(X) > 0$. Then there exists a Radon probability measure $\mu$ and a constant $C$ such that $\mu(X) > 0$ and such that
\begin{equation}\label{BallCondition}
\mu(B) \leq  C r^s \text{ for all closed balls $B$ of radius $r$}.
\end{equation}
Conversely, if $X$ is a compact metric space supporting such a measure $\mu$, then we have $\mathcal{H}^s(X) > 0$.
\end{mylem}
Technically, the converse statement does not appear in Theorem 8.17, but is easily shown to follow from the definition of Hausdorff dimension and a simple calculation similar to the Euclidean version appearing in Theorem 2.7 from Mattila \cite{Mattila15}.
The equation \eqref{BallCondition} is called the $s$-dimensional ball condition. On $\mathbb{F}_q^{\infty}$, the $s$-dimensional ball condition is related to the finiteness of the $s$-energy of $\mu$. The following lemma appears in \cite{Papadimitropoulos10C}:
\begin{mylem}
Let $\mu$ be a Borel probability measure on $\mathbb{F}_q^{\infty}$ satisfying the $s$-dimensional ball condition \eqref{BallCondition}. If $t < s$, then the $t$-energy
\begin{equation}\label{tEnergy}
\iint |x - y|^{-t} \, d \mu(x) \, d \mu(y)
\end{equation}
is finite for any $t < s$.

Conversely, if the $t$-energy \eqref{tEnergy} is finite, then there exists a set $A \subset \mathbb{F}_q^{\infty}$ such that $\mu(A) > 0$ and such that the restriction $\mu|_A$ satisfies the $t$-dimensional ball condition.
\end{mylem}
There is also a Fourier-analytic expression for the $t$-energy. This lemma can also be found in \cite{Papadimitropoulos10C}.
\begin{mylem}
If $\mu$ is a probability measure on $\mathbb{F}_q^{\infty}$,
\[\iint |x - y|^{-t} d \mu(x) \, d \mu(y) = \frac{1 - q^t}{1 - q^{t-1}} \int |\hat \mu(\xi)|^2 |\xi|^{t - 1} \, d \xi.\]
Therefore, the $t$-energy of $\mu$ is finite if and only if $\int |\hat \mu(\xi)|^2 |\xi|^{t - 1} \, d \xi$ is finite.
\end{mylem}
We are now ready to define the Fourier dimension of a compact subset $E \subset \mathbb{F}_q^{\infty}$.
\begin{mydef}
The Fourier dimension of $E$ is the supremum over all real numbers $s$ such that there exists a measure $\mu_s$ supported on $E$ such that
\begin{equation}\label{FourierDimension}
|\hat \mu_s(\xi)| \leq C_s |\xi|^{-s/2}.
\end{equation}
\end{mydef}
It is easy to see that any measure satisfying \eqref{FourierDimension} will have finite $s$-energy- thus a set of Fourier dimension $s_0$ will support a measure with finite $s$-energy for any $s < s_0$.

Combining all of these facts gives the following simple statement:
\begin{mylem}\label{FourierImpliesBall}
Suppose $\mu$ is a measure supported on a compact set $E \subset \mathbb{F}_q^{\infty}$ such that $|\hat \mu(\xi)| \leq C |\xi|^{\beta/2}$ for some constant $C$ and all $\xi \in \widehat{\mathbb{F}_q^{\infty}}$. Then there exists a set $A$ such that $\mu(A) > 0 $ and such that $\mu|_A$ satisfies the $\alpha$-dimensional ball condition for any $\alpha < \beta$. 
\end{mylem}
\section{A Varnavides-type theorem for thin subsets of $\mathbb{F}_q^d$}
Varnavides's theorem \cite[Theorem 10.9]{TaoVu06} gives a quantitative statement about the number of three-term arithmetic progressions in large subsets of $\{1,\ldots,N\}$. We will prove a similar result for the finite group $\mathbb{F}_q^d$. This version of Varnavides's lemma was established by Pohoata and Roche-Newton \cite{PohoataRocheNewton19} using the result of Ellenberg and Gijswijt \cite{EllenbergGijswijt17} and the triangle removal lemma. We present a different proof using a simple counting argument instead of the triangle removal lemma. The proof is similar to the standard proof of Varnavides's theorem, and in particular uses the strategy of intersecting with random planes described by Tao and Vu \cite[Exercise 10.1.9]{TaoVu06} to arrive at a quantitative statement for thin sets.
\begin{mylem}[Varnavides's theorem for $\mathbb{F}_q^d$]\label{Varnavides}
For any odd prime $q$, there exists a positive real number $\alpha^*(q) < 1$ with the following property. Let $A \subset \mathbb{F}_q^d$ have at least $q^{\alpha d}$ elements, where $\alpha > \alpha^*(q)$. Then $A$ contains at least $q^{2 - O_q(1 - \alpha))d}$ three-term arithmetic progressions. 
\end{mylem}
\begin{proof}
Let $d$ be a large integer and let $d' < d$ be a parameter that will be chosen later. We will consider planes in $\mathbb{F}_q^d$ of dimension $d'$. Note that because $A$ has at least $q^{\alpha d}$ elements, it follows that an average plane of $\mathbb{F}_q^d$ of dimension $d'$ will contain $q^{\alpha d + d' - d}$ elements of $A$.

The result of Ellenberg and Gijswijt \cite{EllenbergGijswijt17} implies that if $d'$ is sufficiently large, any subset of a ${d'}$-dimensional plane consisting of at least $q^{\alpha_0 d}$ elements will contain a three-term arithmetic progression, where $\alpha_0 < 1$ is a real number depending only on $q$. We will thus look for planes that contain at least $q^{\alpha_0 d}$ elements of $A$.

Let $W$ be the fraction of planes of $\mathbb{F}_q^d$ of dimension $d'$ that contain at most $q^{\alpha_0 d}$ elements of $A$. We will apply a pigeonhole-principle argument in order to obtain an upper bound for $W$.

As discussed above, the average number of elements of $A$ contained in a plane of dimension $d'$ is at least $q^{\alpha d} q^{d' - d}$. On the other hand, the average number of elements of $A$ in such a plane is bounded above by $W q^{\alpha_0 d} + (1 - W) q^{d'}$. This gives the inequality
\[q^{\alpha d + d' - d} \leq W q^{\alpha_0 d'} + (1 - W) q^{d'}.\]
When we isolate $W$ in this inequality, we arrive at the inequality
\begin{equation}\label{WInequality}
W \leq \frac{1 - q^{(\alpha - 1)d}}{1 - q^{(\alpha_0 - 1)d'}}.
\end{equation}
We choose
\begin{equation}\label{dprimeChoice}
d' = \left\lfloor \frac{10000d}{9999} \cdot \frac{1 - \alpha}{1 - \alpha_0} \right\rfloor.
\end{equation}
If $\alpha$ is sufficiently close to $1$ depending on $q$, then we will have $d' < d$. From \eqref{WInequality} and \eqref{dprimeChoice}, we get
\[W \leq \frac{1 - q^{(\alpha - 1)d}}{1 - q^{\frac{10000d}{9999} (\alpha - 1) + 1}}\]
If $d$ is large enough depending on $\alpha$, then $1 - W$ will be larger than $\frac{q^{(\alpha - 1)d}}{2}$ as can be seen by using e.g. the linearization of the function $\frac{1 - x}{1 - y}$ near $(0,0)$. This gives a lower bound on the fraction of $d'$-dimensional planes that contain at least one three-term arithmetic progression of $A$.

We will now obtain an upper bound for the number of planes of dimension $d'$ that contain a specific three-term arithmetic progression $\{x, x + a, x + 2a\}$. Any such plane is a translation of a subspace containing $a$. Furthermore, given any such subspace, a $q^{d' - d}$ fraction of translations of it will contain the three-term arithmetic progression $\{x, x+a, x+ 2a\}$.

It remains to compute the fraction of $d'$-dimensional subspaces of $\mathbb{F}_q^d$ that contain $a$. We will first compute the number of such subspaces; we will then divide by the total number of subspaces of dimension $d'$ in $\mathbb{F}_q^d$. The number of subspaces of dimension $d'$ containing $a$ is $L_1 / L_2$, where $L_1$ is the number of linearly independent collections of $d'$ vectors containing $a$, and $L_2$ is the number of bases of any vector space of dimension $d'$ containing $a$. We will estimate both $L_1$ and $L_2$.

We compute $L_1$ in the following way. Choose $d'-1$ elements of $\mathbb{F}_q^d$ with possible repetitions.  The probability that the first such element $x_1$ lies in the span of $a$ is $q^{1 - d}$. If the first element does not lie in the span of $a$, then the probability that $x_2$ does not lie in the span of $a$ and $x_1$ is $q^{2 - d}$. Continuing in this manner, we see that the probability that there is a linear dependence among $\{a, x_1, \ldots, x_{d-1}\}$ is
\[q^{1 - d} + (1 - q^{1 -d}) q^{2 - d} + \cdots + (1 - q^{1-d})(1 - q^{2 - d}) \cdots (1 - q^{d' - 2 - d})q^{d' - 1 -d}.\]
This is crudely bounded above by $2 \cdot q^{d' - 1 - d}$ by estimating each $(1 - q^{j-d})$ by $1$ and using the geometric series formula. With the choice of $d'$ given in \eqref{dprimeChoice}, this probability approaches $0$ as $d \to \infty$. Thus $L_1 = (1 + o(1)) q^{d (d'-1)}$.

The number of bases containing $a$ of a $d'$-dimensional vector space over $\mathbb{F}_q$ is computed in a similar way. Again, choosing a $d'-1$-element list of elements of $\mathbb{F}_q^{d'}$ uniformly, we get that the probability of a linear dependence is 
\[q^{1 - d'} + (1 - q^{1 - d'})q^{2 - d'} + \cdots + (1 - q^{1 - d'})(1 - q^{2 - d'}) \cdots (1 - q^{-2}) q^{-1}.\]
For any odd prime $q$, this is bounded above by $\frac{1}{2}$, as can be seen again from the geometric series formula. Therefore, the number $L_2$ of bases containing $a$ of a fixed $d'$-dimensional subspace is bounded between $\frac{1}{2} q^{d' (d'-1)}$ and $q^{d' (d' - 1)}$. Combining the estimates for $L_1$ and $L_2$, we see that the number of $d'$-dimensional subspaces of $\mathbb{F}_q^d$ containing $a$ is at most $(2 + o(1)) q^{(d - d')(d' - 1)}$. A similar argument shows that the total number of $d'$-dimensional subspaces of $\mathbb{F}_q^d$ is at least $(1 + o(1)) q^{(d - d') d'}$. By dividing, we discover that the fraction of $d'$-dimensional subspaces containing $a$ is at most $(2 + o(1)) q^{d - d'}$. 

Therefore, the fraction of planes of dimension $d'$ contained in $\mathbb{F}_q^d$ that contain $\{x, x+a, x+ 2 a\}$ is no more than a constant times $q^{2(d' - d)}$. Note that the choice of $d'$ guarantees that $d' = O(1 - \alpha) d$. If we divide the fraction of planes that contain at least one three-term AP in $A$ by the fraction that contain a specific three-term AP, we observe (by absorbing the multiplicative constants into the error term) that there are at least $q^{(2 - O(1 - \alpha))d}$ three-term arithmetic progressions contained in $A$. 
\end{proof}
\section{Arithmetic Progression in Subsets of $\mathbb{F}_q^{\infty}$}
\subsection{Finding approximate arithmetic progressions in $E$}
The idea of the proof is the following: we can find a large number of arithmetic progressions in approximations to $E$ by using the Hausdorff dimension assumption on $E$; then, we can use the Fourier regularity of the measure $\mu$ together with the compactness of the set $E$ in order to locate three-term arithmetic progressions in the set $E$.

Let $\mu$ be a measure satisfying assumption 1 of Theorem \ref{mainthm}. Then we have a measure $\mu'$ obtained by restricting the measure $\mu$ to some subset $E'$ of $E$ that satisfies the ball condition of dimension $\alpha$. 

Let $\mu'_d$ be the pushforward of $\mu$ under the projection $\pi_d$. Suppose $\alpha > \alpha^*(q)$, where $\alpha^*(q)$ is defined as in Lemma \ref{Varnavides}. The ball condition implies that $\mu_d'$ is bounded above by $C_1 q^{-\alpha d}$. Let $K = \mu'(\mathbb{F}_q^{\infty})$. If we then define
\[\mu_d''(x) = \begin{cases} 
\mu_d'(x) & \text{if $\mu_d'(x) > K q^{-d}/2$} \\
   0      & \text{if $\mu_d'(x) \leq K q^{-d}/2$}
\end{cases}\]
A simple pigeonholing argument shows $\sum_{x \in \mathbb{F}_q^d} \mu_d''(x) \geq K/2$. 

We will define our set $A$ to be the support of $\mu_d''$ in $\mathbb{F}_q^d$. Because $\mu_d''(x) \leq C_1 q^{-\alpha d}$ for $x \in A$ and $\sum_{x \in A} \mu_d''(x) \geq K/2$, we have the lower bound $|A| \geq \frac{K q^{\alpha d}}{2 C_1}$. We can absorb the constants by replacing $\alpha$ by something slightly smaller: letting $\gamma$ be a real number such that $\alpha^*(q) < \gamma < \alpha$, we have for sufficiently large $d$ (depending on $\gamma, q, C_1,$ and $K$) that $\mu_d''(x) \leq q^{- \gamma d}$ and $|A| \geq q^{\gamma d}$. 

We apply Lemma \ref{Varnavides} to the set $A$. This lemma guarantees that there are at least $q^{(2 - O_q(1 -\gamma ))d }$ three-term arithmetic progressions contained in $A$. We will define an auxiliary function $g(x) = \sum_{a \in \mathbb{F}_q^{d} \setminus \{0\}} \mu_d(x) \mu_d(x + a) \mu_d(x + 2 a)$. Because there are at least $q^{(2 - O_q(1 - \gamma))d}$ pairs $(x, a)$ such that $\{x, x + a, x + 2a\}$ is contained in $A$, and $\mu_d(x) \geq \frac{K}{2}q^{-d}$ on $A$, we have, by absorbing the constant $K/2$ into the $q^{O_q(1 - \gamma)d}$ term, that 
\begin{equation}\label{gsum}
\sum_{x \in \mathbb{F}_q^d} g(x) \geq q^{(-1 - O_q(1 - \gamma))d}.
\end{equation}
\subsection{Refining the approximate arithmetic progressions}
At this point, we will pause to consider what this says about the set $E$. Recall that $\mu'$ is supported on $E'$, and therefore if $x \in A$, then there exists $z \in \pi_d^{-1}(A)$ in the support of $\mu'$. We have therefore found many triples of points $\{z, z + a_1, z + 2 a_2\}$ where $\pi_d(a_2) = \pi_d(a_1)$. These are triples of $q^{-d}$-separated points that lie within $q^{-d}$ of the elements of a three-term arithmetic progression.

Let $d^* > d$. We will use the Fourier decay assumption on $\mu$ in order to find three-term arithmetic progressions in the support of $\mu$ (and hence in $E$) that are $q^{-d}$ separated, but that lie within $q^{-d^*}$ of a three-term arithmetic progression. Because we can do this for any $d^*$, the compactness of $E$ will guarantee that, as $d^* \to \infty$, some subsequence of these triples will converge to a three-term arithmetic progression, and the $q^{-d}$-separation of the points in each triple of the subsequence will guarantee that this arithmetic progression is nontrivial.

We begin by defining a function $g_{d^*}$ on $\mathbb{F}_q^{d^*}$ in the following way:

\[g_{d^*}(x) := \sum_{\substack{a \in \mathbb{F}_q^{d^*} \\ |a| \geq q^{-d}}} \mu_d(x) \mu_d(x + a) \mu_d(x + 2a).\]

The function $g_{d^*}(x)$ is nonzero when $x, x + a$, and $x + 2a$ lie in the support of $\mu_{d^*}$ for some $|a| \geq q^{-d}$. Thus if we can show that $\sum_{x \in \mathbb{F}_q^{d^*}} g_{d^*}(x)$ is positive, this will show that there exist $q^{-d}$-separated three-term arithmetic progressions in the support of $\mu_{d^*}$, which in turn would establish that there are $q^{-d}$-separated triples in $E$ that lie within $q^{-d^*}$ of a three-term arithmetic progression.

In order to estimate $\sum_{x \in \mathbb{F}_q^{d^*}} g_{d^*}(x)$, we observe that this sum is equal to $\widehat{g_{d^*}}(0)$. Then
\begin{IEEEeqnarray*}{rCl}
\widehat{g_{d^*}}(0) & = & \sum_{\substack{a \in \mathbb{F}_q^{d^*} \\ |a| \geq q^{-d}}} (\mu_{d^*}(\cdot) \mu_{d^*}(\cdot + a) \mu_{d^*}(\cdot + 2a))^{\wedge}(0)  \\
& = & \frac{1}{q^{d^*}} \sum_{\substack{a \in \mathbb{F}_q^{d^*} \\ |a| \geq q^{-d}}} \sum_{\xi_1 \in \widehat{\mathbb{F}_q^{d^*}}} (\mu_{d^*}(\cdot) \mu_{d^*} (\cdot + a))^{\wedge}(- \xi_1) \widehat{\mu_{d^*}}(\xi_1) e_q(\xi_1 \cdot 2a)\\
& = & \frac{1}{q^{2d^*}} \sum_{\substack{a \in \mathbb{F}_q^{d^*} \\ |a| \geq q^{-d}}} \sum_{\xi_1 \in \widehat{\mathbb{F}_q^{d^*}}} \sum_{\xi_2 \in \widehat{\mathbb{F}_q^{d^*}}} \widehat{\mu_{d^*}}(-\xi_1 - \xi_2) \widehat{\mu_{d^*}}(\xi_1) \widehat{\mu_{d^*}}(\xi_2) e_q((2\xi_1 + \xi_2) \cdot a).
\end{IEEEeqnarray*}
Here, $\xi_1$ and $\xi_2$ are elements of $\mathbb{F}_q^{d^*}$. We will write 
\begin{IEEEeqnarray*}{rCl}
\xi_1 & = & (\xi_1^{(1)}, \ldots, \xi_1^{(d^*)}) \\
\xi_2 & = & (\xi_2^{(1)}, \ldots, \xi_2^{(d^*)}).
\end{IEEEeqnarray*}
We will apply the order $d$ decomposition to $a = a' + a''$, $\xi_1 = \xi_1' + \xi_1''$, $\xi_2 = \xi_2' + \xi_2''$, and observe that the condition $|a| \geq q^{-d}$ is equivalent to the statement that $a' \neq 0$. So we can rewrite this sum as 
\begin{IEEEeqnarray*}{CCl}
= & \frac{1}{q^{2d^*}} \sum_{a' \neq 0} \sum_{a''} \sum_{\xi_1', \xi_2'} \sum_{\xi_1'', \xi_2''} & \bigl(\widehat{\mu_{d^*}}(-\xi_1' - \xi_1'' - \xi_2' - \xi_2'') \widehat{\mu_{d^*}}(\xi_1' + \xi_1'') \widehat{\mu_{d^*}}(\xi_2' + \xi_2'') \cdot \\ 
& & \cdot e_q((2 \xi_1' + \xi_2') \cdot a') e_q((2 \xi_1'' + \xi_2'') \cdot a'')\bigr)
\end{IEEEeqnarray*}
We rearrange this sum so that the sums in $a'$ and $a''$ are inside:
\begin{IEEEeqnarray*}{CCl}
= & \frac{1}{q^{2d^*}} \sum_{\xi_1', \xi_2'} \sum_{\xi_1'', \xi_2''} & \widehat{\mu_{d^*}}(-\xi_1' - \xi_1'' - \xi_2' - \xi_2'') \widehat{\mu_{d^*}}(\xi_1' + \xi_1'') \widehat{\mu_{d^*}}(\xi_2' + \xi_2'') \cdot \\
& & \cdot \left( \sum_{a' \neq 0} e_q((2 \xi_1' + \xi_2') \cdot a')\right) \left(\sum_{a''} e_q((2 \xi_1'' + \xi_2'') \cdot a'')\right)
\end{IEEEeqnarray*}
We will start with the sum
\[\sum_{a''} e_q((2 \xi_1'' + \xi_2'') \cdot a'').\]
This sum vanishes if $2 \xi_1'' + \xi_2''$ is nonzero. If $2 \xi_1'' + \xi_2''$ is equal to zero, then each summand is equal to $1$, so the sum is equal to $q^{d^* - d}$, the number of summands. Therefore, we have
\begin{IEEEeqnarray*}{rCCl}
\widehat{g_{d^*}}(0) & = & q^{-d^* - d} \sum_{\xi_1', \xi_2'} \sum_{\xi_1''}& \widehat{\mu_{d^*}}(- \xi_1' - \xi_2' + \xi_1'') \widehat{\mu_{d^*}}(\xi_1' + \xi_1'') \widehat{\mu_{d^*}}(\xi_2' - 2 \xi_1'') \cdot \\
& & & \cdot \left(\sum_{a' \neq 0} e_q((2 \xi_1' + \xi_2') \cdot a') \right)
\end{IEEEeqnarray*}
We will write
\[\widehat{g_{d^*}}(0) = S_0 + S_{\neq 0},\]
where
\[S_0 := q^{-d^* - d} \sum_{\xi_1', \xi_2'} \widehat{\mu_{d^*}}(-\xi_1' - \xi_2') \widehat{\mu_{d^*}}(\xi_1') \widehat{\mu_{d^*}}(\xi_2')\left(\sum_{a' \neq 0} e_q((2 \xi_1' + \xi_2') \cdot a') \right) \]
and
\begin{IEEEeqnarray*}{rCCl}
S_{\neq 0 } & := & q^{-d^* - d} \sum_{\xi_1', \xi_2'} \sum_{\xi_1'' \neq 0} & \widehat{\mu_{d^*}}(-\xi_1' - \xi_2' + \xi_1'') \widehat{\mu_{d^*}}(\xi_1' + \xi_1'') \widehat{\mu_{d^*}}(\xi_2' - 2 \xi_1'') \cdot \\
& & & \cdot \left(\sum_{a' \neq 0} e_q((2 \xi_1' + \xi_2') \cdot a') \right).
\end{IEEEeqnarray*}
\subsection{Estimating the main term}
We will show that $S_0$ will be the main term and $S_{\neq 0}$ is an error term. We will consider $S_0$ first. To do this, we notice that a calculation similar to the one for $\widehat{g_{d^*}}(0)$ allows us to conclude
\[\hat g(0) = q^{-2d} \sum_{\eta_1, \eta_2 \in \widehat{\mathbb{F}_q^d}} \widehat{\mu_d}(-\eta_1 - \eta_2) \widehat{\mu_d}(\eta_1) \widehat{\mu_d}(\eta_2) \sum_{\substack{b \in \mathbb{F}_q^d \\ b \neq 0}} e_q((2 \eta_1 + \eta_2) \cdot b). \]
Because $\xi_1'$ and $\xi_2'$ have absolute value at most $q^d$, it follows that $\xi_1' \cdot x = \xi_1' \cdot \pi_d(x)$ for any $x \in \mathbb{F}_q^{d^*}$. Thus, if $\eta_1$ is the vector $(\xi_1^{(1)}, \ldots, \xi_1^{(d)})$ and $\eta_2$ is the vector $(\xi_2^{(1)}, \ldots, \xi_2^{(d)})$, then $\widehat \mu_{d^*}(\xi_1')$ is equal to $\widehat \mu_d(\eta_1)$, and similarly for $\xi_2'$ and $-\xi_1' - \xi_2'$. Re-indexing the sum in $\xi_1'$ and $\xi_2'$ in $S_0$ by $\eta_1$ and $\eta_2$, and $b = \pi_d(a)$, we observe
\begin{IEEEeqnarray*}{rCl}
S_0 & = & q^{-d^* - d} \sum_{\eta_1, \eta_2} \widehat{\mu_{d}}(- \eta_1 - \eta_2) \widehat{\mu_{d}}(\eta_1) \widehat{\mu_{d}}(\eta_2) \sum_{\substack{b \in \mathbb{F}_q^d \\ b \neq 0}} e_q ((2 \eta_1 + \eta_2) \cdot b) \\
& = & q^{d- d^*} \hat g(0).
\end{IEEEeqnarray*}
We can bound this from below by using \eqref{gsum} to conclude that if $d$ is sufficiently large depending on $\gamma$, $q$, $C_1$ and $K$, then
\begin{equation}\label{S0bound}
S_0 \geq  q^{-d^* - O_q(1 - \gamma)d}.
\end{equation}
\subsection{Estimating the error term}
Now, we must estimate $S_{\neq 0}$. Because $\xi_1''$ is nonzero in this sum, we have in fact that $|\xi_1''| \geq q^{d+1}$. In particular, this means that $|\xi_1''| > \max(|\xi_1'|, |\xi_2'|)$, and thus $|\xi_1' + \xi_1''|$, $|\xi_2' - 2 \xi_1''|$, and $|-\xi_1' - \xi_2' + \xi_1''|$ are all equal to $|\xi_1''|$ by the ultrametric inequality. Applying statement 2 of Theorem \ref{mainthm} and the triangle inequality, we therefore have the estimate
\begin{IEEEeqnarray*}{rCl}
|S_{\neq 0}| & \leq & q^{-d^* - d}\sum_{\xi_1', \xi_2'} \sum_{\xi_1'' \neq 0} |\xi_1''|^{-3 \beta/2} \left| \sum_{a' \neq 0} e_q((2 \xi_1' + \xi_2') \cdot a') \right| \\
             &   =  & q^{-d^* - d}\left(\sum_{\xi_1'' \neq 0} |\xi_1''|^{-3 \beta/2} \right) \left( \sum_{\xi_1', \xi_2'} \left| \sum_{a' \neq 0} e_q((2 \xi_1' + \xi_2') \cdot a') \right| \right).
\end{IEEEeqnarray*}
First, we will estimate
\[\sum_{\xi_1', \xi_2'} \left| \sum_{a' \neq 0} e_q((2 \xi_1' + \xi_2') \cdot a') \right|.\]
The sum 
\[\sum_{a' \neq 0} e_q((2 \xi_1' + \xi_2') \cdot a')\]
will take the value $-1$ if $2 \xi_1' + \xi_2'$ is nonzero, and will take the value $q^d - 1$ otherwise. For a fixed $\xi_1'$, there is exactly one choice of $\xi_2'$ (namely, $-2 \xi_1'$) such that $2 \xi_1' + \xi_2' = 0$. Thus for each $\xi_1'$, we have
\[\sum_{\xi_2'} \left| \sum_{a' \neq 0} e_q((2 \xi_1' + \xi_2') \cdot a')\right| = 2 (q^d - 1),\]
and thus
\begin{equation}\label{aprimesum}
\sum_{\xi_1', \xi_2'} \left| \sum_{a' \neq 0} e_q((2 \xi_1' + \xi_2') \cdot a') \right| = 2q^{d} (q^{d} - 1) \leq 2q^{2d}.
\end{equation}
We will now estimate
\[\sum_{\xi_1'' \neq 0} |\xi_1''|^{-3 \beta/2}.\]
This sum can be rewritten
\[\sum_{j=d+1}^{d^*} (\#\{\xi_1'' : |\xi_1''| = q^j\}) q^{-3 \beta j/2}.\]
Note that $|\xi_1''| = q^j$ whenever $\xi_1''$ has the form 
\[\xi_1'' = (0, \ldots, 0, \xi_1^{(d+1)}, \xi_1^{(d+2)}, \ldots, \xi_1^{(j)}, 0, \ldots, 0)\] 
where $\xi_1^{(j)} \neq 0$. There are $q$ choices for each of $\xi_1^{(d+1)}, \ldots, \xi_1^{(j-1)}$ and $q-1$ choices for $\xi_1^{(j)}$ and thus there are $(q-1)q^{j - d - 1} \leq q^{j - d}$ values of $\xi_1''$ such that $|\xi_1''| = q^j$. Thus
\begin{IEEEeqnarray*}{rCl}
\sum_{\xi_1'' \neq 0} |\xi_1''|^{-3 \beta/2} & \leq & \sum_{j=d+1}^{d^*} q^{j - d} q^{-3 \beta j/2} \\
& =    & q^{-d} \sum_{j=d + 1}^{d^*} q^{j(1 - 3 \beta/2)}\\
& \leq & q^{-d} \sum_{j=d + 1}^{\infty} q^{j(1 - 3 \beta/2)}.
\end{IEEEeqnarray*}
The sum $\sum_{j=d + 1}^{\infty} q^{j(1 - 3 \beta/2)}$ is convergent because of the assumption that $\beta > 2/3$. The geometric series formula gives the estimate $\sum_{j= d+ 1}^{\infty} q^{j(1 - 3 \beta/2)} \leq C_{q, \beta} q^{d(1 - 3 \beta/2)}$. Thus we get
\begin{equation}\label{xi1doubleprimesum}
\sum_{\xi_1'' \neq 0} |\xi_1''|^{-3 \beta/2} \leq C_{q, \beta} q^{-3 d \beta/2}.
\end{equation}
Combining \eqref{aprimesum} and \eqref{xi1doubleprimesum} and absorbing the constant $2$ into $C_{q, \beta}$, we get
\begin{equation}\label{Sneq0bound}
S_{\neq 0 } \leq C_{q, \beta} q^{- d^* - d - 3 d\beta/2 + 2d} = C_{q, \beta} q^{-d^* + (1 - 3 \beta/2)d}.
\end{equation}

Combining \eqref{S0bound} and \eqref{Sneq0bound}, we arrive at the estimate
\[\widehat{g_{d^*}}(0) \geq q^{- d^* - O_q(1 - \gamma) d} - C_{q,\beta} q^{-d^* + (1 - 3 \beta/2) d}\]
for $d$ sufficiently large depending on $C_1, q, K,$ and $\gamma$. But recall that $\gamma$ can be chosen arbitrarily close to $\alpha$. So if $\alpha$  and $\gamma$ are sufficiently large that the $O_q(1 -\gamma)$ term is less than $3 \beta/2 - 1$, then $\widehat{g_{d^*}}(0)$ will be positive provided $d$ is large enough, completing the proof under conditions 1 and 2.

Notice that if $\mu$ satisfies assumption 2 of Theorem \ref{mainthm}, then Lemma \ref{FourierImpliesBall} implies that assumption 1 holds with $\alpha$ arbitrarily close to  $\beta$ and with some value of $C_1(\alpha)$.  If $\beta$ is sufficiently close to $1$, then choosing $\alpha$ and $\gamma$ sufficiently close to $\beta$ will guarantee that $O_q(1 - \gamma)$ term will be less than $3 \beta/2 - 1$, so the proof is complete in this case as well. \qedhere
\section{Concluding Remarks}
We crucially used the result of Ellenberg-Gijswijt \cite{EllenbergGijswijt17} in the proof of Lemma \ref{Varnavides}. For this reason, the proof described here does not apply to the Euclidean setting. As stated before, Shmerkin \cite{Shmerkin17} has provided a counterexample to Theorem \ref{mainthm} in $\mathbb{R}$. 

The only use of the Fourier decay condition occurred in the estimate of the term $S_{\neq 0}$. 
\bibliographystyle{myplain}
\bibliography{Fourier_Dimension_3_Term_APs}
\end{document}